\ifx\shlhetal\undefinedcontrolsequence\let\shlhetal\relax\fi

\documentclass[11pt]{amsart}

\usepackage{amsmath}
\usepackage{amssymb}

\newtheorem{theorem}{Theorem}[section]
\newtheorem{claim}[theorem]{Claim}

\newtheorem{lemma}[theorem]{Lemma}

\newtheorem{proposition}[theorem]{Proposition}

\theoremstyle{definition}
\newtheorem{definition}[theorem]{Definition}

\newtheorem{question}[theorem]{Question}

\theoremstyle{remark}
\newtheorem{remark}[theorem]{Remark}
\newtheorem{discussion}[theorem]{Discussion}

\newcount\skewfactor
\def\mathunderaccent#1#2 {\let\theaccent#1\skewfactor#2
\mathpalette\putaccentunder}
\def\putaccentunder#1#2{\oalign{$#1#2$\crcr\hidewidth
\vbox to.2ex{\hbox{$#1\skew\skewfactor\theaccent{}$}\vss}\hidewidth}}

\def\smallbox#1{\leavevmode\thinspace\hbox{\vrule\vtop{\vbox
   {\hrule\kern1pt\hbox{\vphantom{\tt/}\thinspace{\tt#1}\thinspace}}
   \kern1pt\hrule}\vrule}\thinspace}

\newcommand{\cf}{{\rm cf}}

\newcommand{\then}{{\underline{then}}}
\newcommand{\Then}{{\underline{Then}}}

\def\qedref#1{$\qed_{\reforiginal{#1}}$}

\title{Many normal measures}
\author{Shimon Garti}
\address{Institute of Mathematics,
 The Hebrew University of Jerusalem,
 Jerusalem 91904, Israel}
\email{shimon.garty@mail.huji.ac.il}


\begin{document}
\let\labeloriginal\label
\let\reforiginal\ref

\begin{abstract}
We characterize the situation of having at least $(2^\kappa)^+$-many normal ultrafilters on a measurable cardinal $\kappa$. We also show that if $\kappa$ is a compact cardinal then $\kappa$ carries $(2^\kappa)^+$-many $\kappa$-complete ultrafilters, each of which extends the club filter on $\kappa$.
\end{abstract}

\maketitle

\newpage

\section{introduction}

The number of normal measures on large cardinals has been investigated extensively in recent years. We know (by Kunen, in \cite{MR0277346}) that it is consistent to have a measurable cardinal which carries just one normal measure on it. 
The opposite direction is taken by Kunen-Paris (see \cite{MR0277381}), who proved the consistency of having $2^{2^\kappa}$-many normal measures on a measurable cardinal $\kappa$. Since then, many papers have dealt with all the values in between (see \cite{MR0344123} under some assumptions on the Mitchell order, \cite{MR2299507} for having $\kappa^+$-many normal measures on the measurable cardinal $\kappa$ without further assumptions, \cite{MR820124}, and recently \cite{MR2548481}, which covers for the measurable cardinal $\kappa$ all the possibilities between $2$ and $\kappa^{++}$ in a uniform manner). However, a measurable cardinal may carry just one normal measure.
On the other hand, if $\kappa$ is supercompact then there are many normal ultrafilters on it (see \cite{MR0347607}). 

We characterize, in this paper, the existence of many normal measures on a measurable cardinal. By ``many" we mean at least $(2^\kappa)^+$, which is an upper bound if the GCH holds.

The idea is pretty simple. By a basic result of Scott (from \cite{MR0166107}), we have a canonical way to normalize a $\kappa$-complete ultrafilter on $\kappa$. If we start with two distinct $\kappa$-complete ultrafilters, we may get the same normal ultrafilter after normalizing both of them. But if we have many different $\kappa$-complete ultrafilters, and some property of the function which represents $\kappa$ in the ultraproduct, then we will be able to show that $\kappa$ carries many different normal ultrafilters. 

We indicate that these results apply to every measurable cardinal, but if one assumes that $\kappa$ is compact, then many $\kappa$-complete ultrafilters on $\kappa$ are available. Moreover, we can show that many of them include the club filter on $\kappa$. This fact does not suffice for normality, since we need an extra property on these ultrafilters, but it is a meaningful step toward proving that there are many normal ultrafilters on every compact cardinal $\kappa$.

Some comments about notation. We shall use $\kappa, \lambda, \mu, \tau$ for infinite cardinals, typically $\tau=(2^\kappa)^+$. We use $\alpha, \beta, \gamma, \delta, \zeta, \eta, \xi$ for ordinals. $D$ will be a filter and $U$ will be an ultrafilter. For an uncountable regular cardinal $\kappa$ we denote the club filter on $\kappa$ by $D_\kappa$. A filter is $\kappa$-complete if it is closed under intersections of \emph{fewer than} $\kappa$ many members from it. A filter is normal if it is closed under diagonal intersections of length $\kappa$. We assume, from now on, that all the filters and ultrafilters are non-principal. For a regular cardinal $\kappa$ and $A,B\subseteq\kappa$, $A$ and $B$ are almost disjoint if $|A\cap B|<\kappa$.

If $\kappa$ is measurable and $U$ is a $\kappa$-complete ultrafilter on $\kappa$, then we can form the ultraproduct ${\rm \bf V}^\kappa/U$. Being well-founded, there is a transitive model ${\rm M}$ which is isomorphic to ${\rm \bf V}^\kappa/U$. We denote the elementary embedding of ${\rm \bf{V}}$ into ${\rm \bf V}^\kappa/U$ by $\jmath$, so the picture is $\jmath : {\rm \bf{V}} \hookrightarrow {\rm \bf {V}}^\kappa / U \cong M$, and we shall not distinguish between ${\rm \bf V}^\kappa/U$ and ${\rm M}$.

Recall that $\kappa$ is the first ordinal moved by $\jmath$, and we call it the critical point. Dealing with many ultrafilters, we denote the elementary embedding which $U$ generates by $\jmath_U$, and if $\delta=\jmath_U(\kappa)$ then we say that $\delta$ is the critical value of $\jmath_U$.

If $U$ is a $\kappa$-complete ultrafilter on $\kappa$, then $\kappa\in M$ and corresponds to an equivalence class in ${\rm \bf V}^\kappa/U$. One can choose a representative for this class, i.e., a function $f:\kappa\rightarrow\kappa$ such that $[f]_U=\kappa$. We shall use the following:

\begin{lemma}
\label{iidd}
$\kappa$-complete ultrafilters and the identity function. \newline 
Let $\kappa$ be a measurable cardinal, $U$ a $\kappa$-complete ultrafilter on $\kappa$. Let $g:\kappa\rightarrow\kappa$ be the identity (i.e., $g(\alpha)=\alpha$ for every $\alpha<\kappa$).
\begin{enumerate}
\item [$(a)$] $[g]_U\geq\kappa$ in ${\rm \bf V}^\kappa/U$.
\item [$(b)$] $U$ is normal iff $[g]_U=\kappa$.
\end{enumerate}
\end{lemma}

\hfill \qedref{iidd}

For a regular cardinal $\kappa$, $\mathcal{A}\subseteq \mathcal{P}(\kappa)$ is almost disjoint if $A,B\in\mathcal{A}\Rightarrow|A\cap B|<\kappa$. We quote the following theorem from \cite{MR597342}, to be used in Theorem \ref{ddkappa} below. The proof appears in \cite{MR597342}, page 48:

\begin{theorem}
\label{kkkunen} If $\kappa\geq\aleph_0$ and $2^{<\kappa}=\kappa$, then there is an almost disjoint family $\mathcal{A}\subseteq\mathcal{P}(\kappa)$ with $|\mathcal{A}|=2^\kappa$.
\end{theorem}

The paper contains two additional sections. In the first one we deal with measurable cardinals, and we give a general characterization for having at least $(2^\kappa)^+$-many normal measures on the measurable cardinal $\kappa$. In the second we deal with compact cardinals, with goal showing the plausibility of having many normal measures on every compact cardinal.

\medskip 

I would like to thank the referee for a careful reading of the manuscript, mathematical corrections and many helpful suggestions which greatly improved the clarity of the paper. I also thank Shani Ben-David for some comments.

\newpage 

\section{Many normal measures on measurable cardinals}

Let us start with the following lemma, which is due to Scott (\cite{MR0166107}). Originally, he used it for showing that if $\kappa$ is an uncountable measurable cardinal (in the sense that $\kappa$ carries a $\kappa$-complete ultrafilter on it), then there is also a normal ultrafilter on $\kappa$.

\begin{proposition}
\label{sscot}
Normalizing an ultrafilter. \newline 
Suppose $\kappa$ is an uncountable measurable cardinal, and $U$ is a $\kappa$-complete ultrafilter on $\kappa$. 
Assume $f:\kappa\rightarrow\kappa$ satisfies $[f]_U = \kappa$. 
\newline 
Set $U^* = \{B\subseteq\kappa:f^{-1}(B)\in U\}$. 
\Then\ $U^*$ is a normal ultrafilter on $\kappa$.
\end{proposition}

\hfill \qedref{sscot}

Starting with two distinct ultrafilters $U_\alpha, U_\beta$, we can normalize and get $U_\alpha^*, U_\beta^*$. But even if $U_\alpha \neq U_\beta$ it may happen that $U_\alpha^* = U_\beta^*$. This setting happens, for instance, in ${\rm \bf L}[U]$ (see \cite{MR0277346}). It is proved there that there exists a measurable cardinal $\kappa$ which carries just one normal measure $U_*$. Take two distinct $\kappa$-complete ultrafilters $U_\alpha$ and $U_\beta$ on $\kappa$, and normalize them. As $U_*$ is the only normal measure on $\kappa$, one must conclude that $U_\alpha\neq U_\beta$ yet $U_\alpha^*=U_\beta^*=U_*$.

Moreover, it is known that in the model ${\rm \bf L}[U]$ there are $\kappa^+$-many nonprincipal $\kappa$-complete ultrafilters on $\kappa$ (see \cite{MR1940513}, Corollary 19.22). It means that the Scott process of normalizing projects all of them onto the same normal ultrafilter.
The key point which enables us to prove the converse is based on the following:

\begin{definition}
\label{ubounded}
Bounded functions. \newline 
Let $U$ be a $\kappa$-complete ultrafilter on $\kappa$, $\theta\leq\kappa, \theta\geq 1$ a (finite or infinite) cardinal, and $f:\kappa\rightarrow\kappa$.
\begin{enumerate}
\item [$(\aleph)$] $f$ is $\theta$-bounded if $|f^{-1}(\beta)|<\theta$ for every $\beta<\kappa$.
\item [$(\beth)$] $f$ is $U$-bounded if $|f(A)\cap f(\kappa\setminus A)| < \kappa$ for every $A\in U$.
\item [$(\gimel)$] $f$ is $(U,U^*)$-bounded if $f(A)\cap f(\kappa\setminus A) \notin U^*$ for every $A\in U$ ($U^*$ is the normalization of $U$ according to Proposition \ref{sscot}).
\item [$(\daleth)$] $f$ is stationarily-bounded if $f^{-1}(\beta)$ is not a stationary set for every $\beta<\kappa$.
\end{enumerate}
\end{definition}

We observe that parts $(\aleph)$ and $(\beth)$ of the definition are incompatible, but part $(\gimel)$ gives a weaker assumption on $f$ than part $(\beth)$. For this, see Remark \ref{rrrr} below.
We shall use part $(\daleth)$ in the next section, when we deal with compact cardinals.
The following is the first characterization that we can prove:

\begin{theorem}
\label{mt}
Normal ultrafilters and $U$-boundedness. \newline 
Let $\kappa$ be a measurable cardinal, $\tau=(2^\kappa)^+$. \newline 
The following conditions are equivalent:
\begin{enumerate}
\item [$(\alpha)$] $\kappa$ carries (at least) $\tau$-many normal ultrafilters.
\item [$(\beta)$] $\kappa$ carries $\tau$-many $\kappa$-complete ultrafilters $(U_\alpha:\alpha<\tau)$, such that for every $\alpha<\tau$ one can choose a $U_\alpha$-bounded $f_\alpha:\kappa\rightarrow\kappa$ with $[f_\alpha]_{U_\alpha}=\kappa$.
\end{enumerate}
Moreover, one can replace the (strong) requirement of $U_\alpha$-boundedness by the (weaker) assumption of $(U_\alpha,U_\alpha^*)$-boundedness.
\end{theorem}

\par \noindent \emph{Proof}. \newline 
$(\alpha)\Rightarrow(\beta)$ is immediate, since every one-to-one function is $U$-bounded (and the identity is one-to-one). We shall prove $(\beta)\Rightarrow(\alpha)$. Let $f_\alpha$ be a $U_\alpha$-bounded function such that $[f_\alpha]_{U_\alpha}=\kappa$, for every $\alpha<\tau$.

Since $\tau=\cf(\tau)>2^\kappa$, there exists $f:\kappa\rightarrow\kappa$ such that $f=f_\alpha$ for $\tau$-many (and without loss of generality, every) $\alpha<\tau$. Use Proposition \ref{sscot} to normalize each $U_\alpha$ and to get a normal ultrafilter $U_\alpha^*$ for every $\alpha<\tau$. We claim now that $\alpha<\beta<\tau\Rightarrow U_\alpha^*\neq U_\beta^*$.

Recall that $U_\alpha\neq U_\beta$, and being ultrafilters one can choose $A\in U_\alpha$ such that $\kappa\setminus A\in U_\beta$. Denote $f(A)$ by $B_0$ and $f(\kappa\setminus A)$ by $B_1$. The assumption on $f_\alpha$ (which equals $f$) implies $|B_0\cap B_1|<\kappa$. Now $B_0\in U_\alpha^*$ (since $f^{-1}(B_0)\supseteq A$) and $B_1\in U_\beta^*$ (since $f^{-1}(B_1)\supseteq\kappa\setminus A$), but $B_1\notin U_\alpha^*$ (since $B_0\in U_\alpha^*$ and $|B_0\cap B_1|<\kappa$). It means that $U_\alpha^*\neq U_\beta^*$, and we are done. The same argument gives also the same conclusion under the weaker assumption of $(U_\alpha,U_\alpha^*)$-boundedness, so the proof is complete.

\hfill \qedref{mt}

\begin{remark}
\label{rrrr}
We defined (in \ref{ubounded} above) the notions of $\theta$-boundedness and $U$-boundedness. These notions are really distinct, as a representing function might be $\theta$-bounded (even for $\theta=3$) but not $U$-bounded (e.g., if there is an $A\in U$ such that $|A|=|\kappa\setminus A|=\kappa$ and for every $\beta<\kappa$ there exist unique $\alpha_0\in A,\alpha_1\in \kappa\setminus A$ so that $f(\alpha_0)=f(\alpha_1)=\beta$). Similarly, a representing function might be $U$-bounded but not $\kappa$-bounded (e.g., if there is some $A$ of cardinality $\kappa$, $A\notin U$, and $f(A)=\{\beta\}$ for some $\beta<\kappa$).

It means that parts $(\aleph)$ and $(\beth)$ in Definition \ref{ubounded} are incompatible. Clearly, part $(\gimel)$ is weaker than $(\beth)$, as $U^*$ is nonprincipal and $\kappa$-complete hence contains no set of size less than $\kappa$.
\end{remark}

Nevertheless, we have another characterization (similar to the previous one), based on $\theta$-bounded functions:

\begin{theorem}
\label{mmt}
Normal ultrafilters and $\theta$-boundedness. \newline 
Let $\kappa$ be a measurable cardinal, $\tau=(2^\kappa)^+$, and $\theta_\alpha<\kappa$ for every $\alpha<\tau$. \newline 
The following conditions are equivalent:
\begin{enumerate}
\item [$(\alpha)$] $\kappa$ carries (at least) $\tau$-many normal ultrafilters.
\item [$(\beta)$] $\kappa$ carries $\tau$-many $\kappa$-complete ultrafilters $(U_\alpha:\alpha<\tau)$, such that for every $\alpha<\tau$ one can choose a $\theta_\alpha^+$-bounded $f_\alpha:\kappa\rightarrow\kappa$ with $[f_\alpha]_{U_\alpha}=\kappa$.
\end{enumerate}
\end{theorem}

\par \noindent \emph{Proof}. \newline 
Since $\tau=\cf(\tau)>\kappa$, the correspondence $\alpha\mapsto\theta_\alpha$ yields an unbounded subset $T\subseteq\tau$ and some $\theta<\kappa$ so that $\alpha\in\ T\Rightarrow\theta_\alpha\equiv\theta$, so without loss of generality $\theta_\alpha\equiv\theta$ for every $\alpha<\tau$.
As in Theorem \ref{mt}, we need only the direction of $(\beta)\Rightarrow(\alpha)$.
Fix, for awhile, an ordinal $\alpha<\tau$.
Denote $f_\alpha^{-1}(\beta)$ by $A_\beta$, for every $\beta<\kappa$. Enumerate the members of $A_\beta$ by $\{a_\gamma^\beta:\gamma<\theta\}$, for each $\beta<\kappa$ (if $|A_\beta|<\theta$ use repetitions, and if $A_\beta=\emptyset$ ignore it). For every $\gamma<\theta$ set $C_\gamma=\{a_\gamma^\beta:\beta<\kappa\}$.

It follows that $\kappa=\bigcup\{C_\gamma:\gamma<\theta\}$. $\kappa\in U_\alpha$ and $U_\alpha$ is $\kappa$-complete, so there is an ordinal $\gamma<\theta$ such that $C_\gamma\in U_\alpha$. By the construction, $f_\alpha\upharpoonright C_\gamma$ is one-to-one. Our purpose is to show that there exists a one-to-one function on $\kappa$ (not just on $C_\gamma$) $h_\alpha\in[f_\alpha]_{U_\alpha}$.

Decompose $C_\gamma$ into two sets, say $C_\gamma=C_0^\gamma\cup C_1^\gamma$, such that $C_0^\gamma\cap C_1^\gamma=\emptyset$ and $|C_0^\gamma|=|C_1^\gamma|=\kappa$. Without loss of generality, $C_0^\gamma\in U_\alpha$ and $C_1^\gamma\notin U_\alpha$. Let $B=\kappa\setminus\{f_\alpha(\xi):\xi\in C_0^\gamma\}$. Notice that $|B|=\kappa$ since $f_\alpha\upharpoonright C_1^\gamma\subseteq B$, and $f_\alpha$ is one-to-one on $C_1^\gamma$.

Choose a one-to-one map $g_\alpha:(\kappa\setminus C_\gamma)\cup C_1^\gamma\rightarrow B$, and define $h_\alpha=f_\alpha\upharpoonright C_0^\gamma\cup g_\alpha$. It follows that $h_\alpha$ is a one-to-one function from $\kappa$ into $\kappa$, and $h_\alpha\in[f_\alpha]_{U_\alpha}$, so $[h_\alpha]_{U_\alpha}=\kappa$.

Repeat this process for every $\alpha<\tau$, and use the fact that $\tau=\cf(\tau)>2^\kappa$ to get a fixed one-to-one function $h$ and $\tau$-many ultrafilters with $h$ as their representing function of $\kappa$. Since $h$ is one-to-one we have $h(A)\cap h(\kappa\setminus A)=\emptyset$ for every $A\subseteq\kappa$. Hence $h$ is $U_\alpha$-bounded for every $\alpha<\tau$. Now employ the previous theorem to conclude that there are (at least) $\tau$-many normal ultrafilters on $\kappa$.

\hfill \qedref{mmt}

\begin{remark}
\label{llambda}
The above theorems focus on the case of $\tau=(2^\kappa)^+$, but $\tau$ can be replaced, verbatim, by any $\lambda\geq\cf(\lambda)\geq\tau$.
\end{remark}

\newpage 

\section{Compact cardinals and the club filter}

Suppose $\kappa$ is a compact cardinal. An old open problem (which appears in \cite{MR1940513} and in \cite{MR1994835}) is whether $\kappa$ carries more than one normal measure. Let us try to explain the background behind this problem.

By a fundamental paper of Magidor (see \cite{MR0429566}), the first compact cardinal has an identity crisis. Magidor proved that on one hand it is consistent that the first compact cardinal is also the first measurable cardinal. On the other hand, it is consistent that the first compact cardinal is the first supercompact cardinal. So from this point of view, the compact cardinal is flexible. But we can judge the compact cardinals from a different point of view.

As noted above, a measurable cardinal may carry just one normal measure. A supercompact cardinal carries many normal measures (see \cite{MR0347607}). 
The problem of how many normal measures exist on a compact cardinal can be viewed as another aspect of the identity crisis. 
Most of the papers which decrease the number of normal measures on a measurable cardinal, mentioned in the introduction, make use of the universe ${\rm \bf L}[U]$ (when $U$ is a normal ultrafilter on $\kappa$). This includes the pivotal result of Kunen (which gives the consistency of having just one normal measure on a measurable cardinal).
We know that there are no compact cardinals in ${\rm \bf L}[U]$ (see \cite{MR0211872}). Yet it is an open problem whether every compact cardinal carries at least two normal measures.

The first claim below shows that if $U$ is a $\kappa$-complete ultrafilter which extends $D_\kappa$, then a simple property of the representing function of $\kappa$ in the ultraproduct suffices for implying the normality of $U$. This claim does not assume compactness. The main theorem of this section is that many ultrafilters extending $D_\kappa$ are available for every compact cardinal.

\begin{claim}
\label{kappabound}
Stationary-boundedness and normality. \newline 
Suppose $\kappa$ is an uncountable measurable cardinal, and $U$ is a $\kappa$-complete ultrafilter on $\kappa$ so that $U \supseteq D_\kappa$. Assume there is a stationarily-bounded function $f:\kappa\rightarrow\kappa$ so that $[f]_U=\kappa$.
\Then\ $U$ is normal. \newline 
Moreover, if we call $\beta$ a good ordinal when $f^{-1}(\beta)$ is not stationary, \then\ the conclusion holds when there is some $\beta_*<\kappa$ so that $\beta\in[\beta_*,\kappa)\Rightarrow\beta$ is a good ordinal and $[f]_U=\kappa$.
\end{claim}

\par \noindent \emph{Proof}. \newline 
Let $g$ be the identity function on $\kappa$.
Define $S_0=\{\alpha<\kappa:f(\alpha)\geq g(\alpha)\}$ and $S_1=\{\alpha<\kappa:f(\alpha)<g(\alpha)\}$. If $S_1\in U$ then $S_1$ is stationary (since $U\supseteq D_\kappa$). By Fodor's lemma there exists a stationary subset $S\subseteq S_1$ such that $f$ is constant on $S$, contrary to the assumptions of the claim.

Hence $S_1\notin U$, so $S_0\in U$. It follows that $g\leq_Uf$. On the other hand, $[f]_U=\kappa$ and $[g]_U\geq\kappa$ (due to Lemma \ref{iidd}(a)), so $f\leq_Ug$. Together, $f\equiv_Ug$, and by Lemma \ref{iidd}(b) we infer that $U$ is normal.

For the additional part of the claim, suppose there is just one bad ordinal $\beta$. It means that $f^{-1}(\beta)$ is stationary, but notice that this set does not belong to $U$ (otherwise, $[f]_U=\beta$). Hence there is $h\in[f]_U$ such that $h$ coincides with $f$ on $\kappa\setminus f^{-1}(\beta)$, and $h\upharpoonright f^{-1}(\beta)$ is one-to-one. Since $h\upharpoonright \kappa\setminus f^{-1}(\beta)\equiv f\upharpoonright \kappa\setminus f^{-1}(\beta)$ we get a stationarily-bounded function. By the $\kappa$-completeness of $U$, this process can be carried out for fewer than $\kappa$-many bad ordinals, so we are done.

\hfill \qedref{kappabound}

Recall that a family $\mathcal{C}$ of subsets of $\kappa$ is $\kappa$-independent if given disjoint subcollections $C_0,C_1\subseteq\mathcal{C}$ each of size less than $\kappa$, $\bigcap\{x:x\in C_0\} \cap \bigcap\{\kappa\setminus x:x\in C_1\}\neq\emptyset$. It is known that if $\kappa$ is inaccessible then there exists a $\kappa$-independent family of size $2^\kappa$ (The proof goes back to Hausdorff, see \cite{MR597342}, p. 288. In order to get $\kappa$-independence one has to replace $\omega$ there by $\kappa$).

A $\kappa$-independent family on $\kappa$ of size $2^\kappa$ yields $2^{2^\kappa}$-many distinct filters. Indeed, let $\mathcal{C}$ be such a family. Let $f$ be a function from $\mathcal{C}$ into $\{0,1\}$. Set $F=\{x:x\in\mathcal{C},f(x)=0\}\cup \{\kappa \setminus x:x\in\mathcal{C},f(x)=1\}$. Insofar as $\mathcal{C}$ is $\kappa$-independent, $F$ is a base for a $\kappa$-complete filter. Clearly, $f\neq g$ create two distinct filters. Inasmuch as there are $2^{2^\kappa}$-many functions from $\mathcal{C}$ into $\{0,1\}$ we have $2^{2^\kappa}$-many distinct $\kappa$-complete filters.

One of the defining properties of compact cardinals is the extension property of $\kappa$-complete filters into $\kappa$-complete ultrafilters. Basically, one can create $2^{2^\kappa}$-many $\kappa$-complete ultrafilters on a compact cardinal in the traditional way. Namely, first create a $\kappa$-independent family of size $2^\kappa$, which yields $2^{2^\kappa}$-many different $\kappa$-complete filters. Extending each filter into a $\kappa$-complete ultrafilter gives the desired result. 

It is not clear if the above method can be employed in order to create many $\kappa$-complete ultrafilters on a compact cardinal $\kappa$, each of which extends the club filter $D_\kappa$. Nonetheless, the following theorem suggests a way to define such ultrafilters:

\begin{theorem}
\label{ddkappa}
Many club ultrafilters on a compact cardinal. \newline 
Suppose $\kappa$ is a compact cardinal, and let $\tau$ be $(2^\kappa)^+$. \newline 
\Then\ there exist $\tau$-many distinct ultrafilters $(U_\alpha:\alpha<\tau)$ such that $D_\kappa\subseteq U_\alpha$ for every $\alpha<\tau$.
\end{theorem}

\par \noindent \emph{Proof}. \newline 
We shall prove that for every ordinal $\delta<(2^\kappa)^+$ there is an ultrafilter $U_\delta$ so that $D_\kappa\subseteq U_\delta$ and $\jmath_{U_\delta}(\kappa)>\delta$. By induction on $\tau$ this yields the desired collection of ultrafilters. Let us start with a disjoint partition of $\kappa$ into $\kappa$-many stationary sets, say $\langle S_\alpha:\alpha<\kappa\rangle$. We employ here the celebrated theorem of Solovay, from \cite{MR0290961}. In fact, Solovay proves much more. He shows that every stationary subset of $\kappa$ can be partitioned into $\kappa$-many disjoint stationary subsets.

Let $\langle t_\beta:\beta<\kappa\rangle$ be an enumeration of $[\kappa]^{<\kappa}$ in such a way that each $t_\beta$ appears stationarily many times (one can use the sequence $\langle S_\alpha:\alpha<\kappa\rangle$ to scatter the members of $[\kappa]^{<\kappa}$ on these stationary sets). Finally, choose an almost disjoint family $\mathcal{A}=\{A_\gamma:\gamma<2^\kappa\}$ in $[\kappa]^\kappa$. Such a family exists, due to Theorem \ref{kkkunen}, upon noticing that $\kappa$ is compact hence $2^{<\kappa}=\kappa$.

For every $\gamma<2^\kappa$ we choose a function $f_\gamma:A_\gamma\rightarrow\kappa\setminus\{0\}$ such that $\beta<\kappa\Rightarrow|f_\gamma^{-1}(\beta)|=\kappa$. For each $A_\gamma\in\mathcal{A}$ we define a function $g_{A_\gamma}:\kappa\rightarrow\kappa$ as follows. For every $\alpha<\kappa$ we check if $t_\alpha\cap A_\gamma$ is a singleton $\{x\}$. If so, then we set $g_{A_\gamma}(\alpha)=f_\gamma(x)$. In all other cases we define $g_{A_\gamma}(\alpha)=0$.

Set $G=\{g_{A_\gamma}:\gamma<2^\kappa\}$. It is easily verified that $A\neq B\Rightarrow g_A\neq g_B$. Indeed, assume $A,B\in\mathcal{A}$ and $A\neq B$. By the fact that $A,B$ are almost disjoint, we can choose some $x\in\kappa$ such that $x\in A$ but $x\notin B$. Set $t=\{x\}$, so $t\in[\kappa]^{<\kappa}$. Pick an ordinal $\alpha<\kappa$ so that $t=t_\alpha=\{x\}$. It follows that $t_\alpha\cap A=\{x\}$, hence $g_A(\alpha)=f_\gamma(x)$ (where $\gamma$ is the index of $A$ in the enumeration of $\mathcal{A}$). Recall that the range of $f_\gamma$ is $\kappa\setminus\{0\}$, and conclude that $g_A(\alpha)>0$.

On the other hand, $t_\alpha\cap B=\emptyset$, as $x\notin B$. In particular, $t_\alpha\cap B$ is not a singleton. By the definition of $g_B$ we have $g_B(\alpha)=0$. Consequently, $g_B(\alpha)\neq g_A(\alpha)$, so the functions $g_B$ and $g_A$ are not the same, i.e., $g_A\neq g_B$.

By the above paragraphs, $|G|=2^\kappa$. Let $\langle g_\varepsilon :\varepsilon\leq\delta\rangle$ be an enumeration of $G$. We emphasize that the enumeration is taken along the ordinal $\delta+1$ (whose cardinality is $2^\kappa$), as we shall need this in the sequel.
For every pair $\alpha<\beta\leq\delta$ we define:

$$
X_{\alpha,\beta}=\{\xi<\kappa:g_\alpha(\xi)<g_\beta(\xi)\}.
$$

Notice that $X_{\alpha,\beta}$ includes a stationary subset of $\kappa$ for every $\alpha<\beta\leq\delta$. Indeed, the functions $g_\alpha,g_\beta$ correspond to some $A,B\in\mathcal{A}$. Choose any $x\in B\setminus A$ (recall that $\mathcal{A}$ is a collection of almost disjoint sets), and we know that $\{x\}\in[\kappa]^{<\kappa}$. By the enumeration of $[\kappa]^{<\kappa}$ there is a stationary $S\subseteq\kappa$ so that $\xi\in S\Rightarrow t_\xi=\{x\}$. Consequently, $\xi\in S\Rightarrow g_\alpha(\xi)=0<f_\varepsilon(x)= g_\beta(\xi)$ (where $B=A_\varepsilon$ in the enumeration of $\mathcal{A}$), so $S\subseteq X_{\alpha,\beta}$.

Moreover, $G$ is endowed with the following property. If $G'\subseteq G, |G'|<\kappa$ and $0<\zeta_A<\kappa$ for every $g_A\in G'$, then there exists a stationary set $S\subseteq\kappa$ such that $g_A(\xi)=\zeta_A$ for every $g_A\in G'$ and every $\xi\in S$. 

To prove this, choose a member $x_A\in A$ for every $g_A\in G'$ such that $g_B \in G' \wedge g_B\neq g_A\Rightarrow x_A\notin B$ and $f_{\gamma_A}(x_A)=\zeta_A$ (where unless otherwise specified, $A= A_{\gamma_A}$ in the enumeration of $G$).
Let us try to justify the existence of such $x_A$-s. The cardinality of the set $\bigcup \{B\cap A:g_B\in G'\wedge g_B\neq g_A\}$ is less than $\kappa$ (as it results as a union of less than $\kappa$-many sets, each of size fewer than $\kappa$, and $\kappa$ is regular).

Since $\zeta_A>0$, the set $f_{\gamma_A}^{-1}(\zeta_A)$ is of size $\kappa$, and included in $A$. So choose $x_A\in A$ such that $x_A\notin \bigcup \{B\cap A:g_B\in G'\wedge g_B\neq g_A\}$ and $x_A\in f_{\gamma_A}^{-1}(\zeta_A)$. It follows that $f_{\gamma_A}(x_A)=\zeta_A$ and $x_A\notin B$ whenever $B\neq A$ and $g_B\in G'$.

With the $x_A$-s at hand, set $s=\{x_A:g_A\in G'\}$. Clearly $s\in[\kappa]^{<\kappa}$, so there exists a stationary set $S$ such that $\xi\in S\Rightarrow s=t_\xi$. It follows that $s\cap A=\{x_A\}$ for every $g_A\in G'$. Now take any ordinal $\xi\in S$ and any $A\in\mathcal{A}$ such that $g_A\in G'$. We have $t_\xi\cap A= s\cap A=\{x_A\}$, hence $g_A(\xi)=f_{\gamma_A}(x_A)=\zeta_A$, and we are done.

Consequently, the intersection of fewer than $\kappa$-many sets of the form $X_{\alpha,\beta}$ includes a stationary set. To prove this, assume $F$ is a collection of $\theta$-many sets of the form $X_{\alpha,\beta}$, for some $\theta<\kappa$. Let $F_{ord}$ be the set of all ordinals mentioned in the index of $X_{\alpha,\beta}$ for some $X_{\alpha,\beta}\in F$. Clearly, $|F_{ord}| \leq\theta<\kappa$.

For every pair $(\alpha,\beta)$ such that $\alpha<\beta$ and $X_{\alpha,\beta}\in F$ choose $\zeta_\alpha,\zeta_\beta$ so that $0<\zeta_\alpha<\zeta_\beta<\kappa$. For this end, enumerate the members of $F_{ord}$ in an increasing order, and choose $\zeta_\beta> {\rm sup} \{\zeta_\alpha:\alpha<\beta\}$ by induction, using the fact that $\theta<\kappa=\cf(\kappa)$.

By the above considerations, we can choose $x_\alpha$ for every $\alpha\in F_{ord}$ such that $x_\alpha\in A_\alpha$ (where $A_\alpha$ is the set in $\mathcal{A}$ which corresponds to $g_\alpha$), $x_\alpha\notin A_\beta$ for every $\beta\in F_{ord}\setminus\{\alpha\}$ and $f_{\varepsilon(\alpha)}(x_\alpha)=\zeta_\alpha$ (where $\varepsilon(\alpha)$ is the index of $A_\alpha$ in the enumeration of $\mathcal{A}$). We denote the set $\{x_\alpha:\alpha \in F_{ord}\}$ by $s$.
Since $s\in[\kappa]^{<\kappa}$, there exists a stationary set $S\subseteq\kappa$ such that $\xi\in S\Rightarrow t_\xi=s$. We claim that $S\subseteq X_{\alpha,\beta}$ for every $X_{\alpha,\beta}\in F$.

For this, pick any ordinal $\xi\in S$, and any set $X_{\alpha,\beta}\in F$. We have chosen $\zeta_\alpha<\zeta_\beta$, and $x_\alpha\in A_\alpha\setminus A_\beta$ as well as $x_\beta\in A_\beta\setminus A_\alpha$ such that $f_{\varepsilon(\alpha)}(x_\alpha)=\zeta_\alpha$ and $f_{\varepsilon(\beta)}(x_\beta)=\zeta_\beta$ (Again, $\varepsilon(\alpha)$ is the index of $A_\alpha$ and $\varepsilon(\beta)$ is the index of $A_\beta$ in the enumeration of $\mathcal{A}$). 
It means that $g_\beta(\xi) = \zeta_\beta$. By the same token, $g_\alpha(\xi)=\zeta_\alpha<\zeta_\beta$, so $\xi\in X_{\alpha,\beta}$ by its very definition. As $\xi\in S$ was arbitrary, we conclude that $S\subseteq X_{\alpha,\beta}$, and as $X_{\alpha,\beta}$ was arbitrary we are done.

It follows from the discussion above that the collection $F=\{X_{\alpha,\beta}:\alpha<\beta\leq\delta\}$ is a $\kappa$-complete filter base. Moreover, $F\cup D_\kappa$ is still a $\kappa$-complete filter base. Let $D$ be the filter generated by $F\cup D_\kappa$, and let $U_\delta$ be any $\kappa$-complete ultrafilter which extends $D$.

By the construction, $\alpha<\beta\leq\delta\Rightarrow X_{\alpha,\beta}\in U_\delta$, so the sequence $\langle g_\alpha:\alpha\leq\delta\rangle$ establishes an increasing sequence of functions according to $<_{U_\delta}$. We conclude that $\jmath_{U_\delta}(\kappa)\geq\delta+1$ as required.

\hfill \qedref{ddkappa}

\begin{discussion}
\label{dddd}
Call $U$ a club ultrafilter if $D_\kappa\subseteq U$. Let $(U_\alpha:\alpha<\tau)$ be a family of club ultrafilters on a compact cardinal $\kappa$, stipulating $\tau=(2^\kappa)^+$. Since the number of functions from $\kappa$ into $\kappa$ is just $2^\kappa$ we may assume, without loss of generality, that the representing function $f_\alpha$ of $\kappa$ in ${\rm \bf V}^\kappa/U_\alpha$ is the same $f$ for every $\alpha<\tau$.

This means that $f^{-1}(\beta)$ is not a club set for every $\beta<\kappa$. 
Indeed, if $f^{-1}(\beta)$ is a club for some ordinal $\beta<\kappa$ then $f^{-1}(\beta)\in D_\kappa\subseteq U_\alpha$ (for every $U_\alpha$). Choose one of the $U_\alpha$-s and conclude that $f$ represents the ordinal $\beta$ in ${\rm \bf V}^\kappa/U_\alpha$, which is impossible as $f$ represents $\kappa$ itself.

Moreover, for every stationary subset $S\subseteq\kappa$, it suffices to find one ordinal $\alpha<\tau$ such that $S\in U_\alpha$ to conclude that $S\neq f^{-1}(\beta)$ for every $\beta<\kappa$. The reason is just the same. If $S= f^{-1}(\beta)$ for some $\beta<\kappa$ and $S\in U_\alpha$ for some $\alpha$, then $f$ represents $\beta$ in ${\rm \bf V}^\kappa/U_\alpha$, contrary to the choice of $f$.

So if we could catch all the stationary subsets of $\kappa$ (in the sense that for every $S$ there is some $\alpha$ such that $S\in U_\alpha$) then we could infer that $f$ is stationarily bounded. As $f$ represents $\kappa$ in $U_\alpha$ for every $\alpha$, this gives the conclusion that each $U_\alpha$ is a normal ultrafilter, by virtue of Claim \ref{kappabound}.
\end{discussion}

\hfill \qedref{dddd}

In light of \ref{dddd}, Theorem \ref{ddkappa} gives rise to the possibility that the existence of $\tau$-many club ultrafilters (for $\tau=(2^\kappa)^+$) is a sufficient condition for having $\tau$-many normal ultrafilters on $\kappa$, with no need to assume compactness. We conclude with the following problem:

\begin{question}
\label{ooo}
Let $\kappa$ be a measurable cardinal, $\tau=(2^\kappa)^+$. \newline 
Assume there are $\tau$-many club ultrafilters on $\kappa$. \newline 
Does this imply that $\kappa$ carries $\tau$-many normal ultrafilters?
\end{question}

\newpage 

\bibliographystyle{amsplain}
\bibliography{arlist}

\providecommand{\bysame}{\leavevmode\hbox to3em{\hrulefill}\thinspace}
\providecommand{\MR}{\relax\ifhmode\unskip\space\fi MR }
\providecommand{\MRhref}[2]{%
  \href{http://www.ams.org/mathscinet-getitem?mr=#1}{#2}
}
\providecommand{\href}[2]{#2}
\begin{thebibliography}{10}

\bibitem{MR2299507}
Arthur~W. Apter, James Cummings, and Joel~David Hamkins, \emph{Large cardinals
  with few measures}, Proc. Amer. Math. Soc. \textbf{135} (2007), no.~7,
  2291--2300 (electronic). \MR{2299507 (2008b:03067)}

\bibitem{MR820124}
Stewart Baldwin, \emph{The {$\triangleleft\,$}-ordering on normal
  ultrafilters}, J. Symbolic Logic \textbf{50} (1985), no.~4, 936--952 (1986).
  \MR{820124 (87d:03124)}

\bibitem{MR2548481}
Sy-David Friedman and Menachem Magidor, \emph{The number of normal measures},
  J. Symbolic Logic \textbf{74} (2009), no.~3, 1069--1080. \MR{2548481}

\bibitem{MR1940513}
Thomas Jech, \emph{Set theory}, Springer Monographs in Mathematics,
  Springer-Verlag, Berlin, 2003, The third millennium edition, revised and
  expanded. \MR{1940513 (2004g:03071)}

\bibitem{MR1994835}
Akihiro Kanamori, \emph{The higher infinite}, second ed., Springer Monographs
  in Mathematics, Springer-Verlag, Berlin, 2003, Large cardinals in set theory
  from their beginnings. \MR{1994835 (2004f:03092)}

\bibitem{MR0166107}
H.~J. Keisler and A.~Tarski, \emph{From accessible to inaccessible cardinals.
  {R}esults holding for all accessible cardinal numbers and the problem of
  their extension to inaccessible ones}, Fund. Math. \textbf{53} (1963/1964),
  225--308. \MR{0166107 (29 \#3385)}

\bibitem{MR0277381}
K.~Kunen and J.~B. Paris, \emph{Boolean extensions and measurable cardinals},
  Ann. Math. Logic \textbf{2} (1970/1971), no.~4, 359--377. \MR{0277381 (43
  \#3114)}

\bibitem{MR0277346}
Kenneth Kunen, \emph{Some applications of iterated ultrapowers in set theory},
  Ann. Math. Logic \textbf{1} (1970), 179--227. \MR{0277346 (43 \#3080)}

\bibitem{MR597342}
\bysame, \emph{Set theory}, Studies in Logic and the Foundations of
  Mathematics, vol. 102, North-Holland Publishing Co., Amsterdam, 1980, An
  introduction to independence proofs. \MR{597342 (82f:03001)}

\bibitem{MR0347607}
M.~Magidor, \emph{There are many normal ultrafiltres corresponding to a
  supercompact cardinal}, Israel J. Math. \textbf{9} (1971), 186--192.
  \MR{0347607 (50 \#110)}

\bibitem{MR0429566}
Menachem Magidor, \emph{How large is the first strongly compact cardinal? or
  {A} study on identity crises}, Ann. Math. Logic \textbf{10} (1976), no.~1,
  33--57. \MR{0429566 (55 \#2578)}

\bibitem{MR0344123}
William~J. Mitchell, \emph{Sets constructible from sequences of ultrafilters},
  J. Symbolic Logic \textbf{39} (1974), 57--66. \MR{0344123 (49 \#8863)}

\bibitem{MR0290961}
Robert~M. Solovay, \emph{Real-valued measurable cardinals}, Axiomatic set
  theory ({P}roc. {S}ympos. {P}ure {M}ath., {V}ol. {XIII}, {P}art {I}, {U}niv.
  {C}alifornia, {L}os {A}ngeles, {C}alif., 1967), Amer. Math. Soc., Providence,
  R.I., 1971, pp.~397--428. \MR{0290961 (45 \#55)}

\bibitem{MR0211872}
P.~Vop{\v{e}}nka and K.~Hrb{\'a}{\v{c}}ek, \emph{On strongly measurable
  cardinals}, Bull. Acad. Polon. Sci. S\'er. Sci. Math. Astronom. Phys.
  \textbf{14} (1966), 587--591. \MR{0211872 (35 \#2747)}

\end{thebibliography}

\end{document}